\documentstyle[11pt,psfig]{amsart}

\def\evol/{{\tt evolver}} \def\illiVert/{{\tt illiVert}}
\def\gv/{{\tt geomview}}
\def\hwm/{halfway-model}

\def\fig#1{Fig.~\ref{fig:#1}} 
 
\def\endf#1#2#3{\caption[#1]{{\small #2}} \label{fig:#3} \end{figure}}
\def\wfig#1#2{\psfig{figure=figs/#1.ps,height=#2in}}
\def\figr#1#2#3#4{\begin{figure}[tbp]
 \centerline{\hbox{\wfig{#1}{#2}}}\endf{#3}{#4}{#1}}
\def\figrp#1#2#3#4{\begin{figure}[p]
 \centerline{\hbox{\wfig{#1}{#2}}}\endf{#3}{#4}{#1}}
\def\figrb#1#2#3#4{\begin{figure}[b!]
 \centerline{\hbox{\wfig{#1}{#2}}}\endf{#3}{#4}{#1}}
\def\figtwo#1#2#3{\begin{figure}[tbp]
 \centerline{ \wfig{#1A}{2} \hfil \wfig{#1B}{2} }
 \endf{#2}{#3}{#1}}
\def\figtwob#1#2#3{\begin{figure}[b!]
 \centerline{ \wfig{#1A}{2} \hfil \wfig{#1B}{2} }
 \endf{#2}{#3}{#1}}
\def\figthree#1#2#3{\begin{figure}[tbp]
 \centerline{ \wfig{#1A}{1.5} \hfil \wfig{#1B}{1.5} \hfil \wfig{#1C}{1.5} }
 \endf{#2}{#3}{#1}}
\begin{document}

\title{``The Optiverse'' and Other Sphere Eversions}

\author{John M. Sullivan}

\address{Mathematics Dept., Univ.~of Illinois, Urbana IL, USA 61801}
\email{jms@@math.uiuc.edu}

\begin{abstract}
For decades, the sphere eversion has been a classic subject for mathematical
visualization.  The 1998 video {\em The Optiverse}
shows geometrically optimal eversions created by minimizing elastic
bending energy.  We contrast these minimax eversions with earlier ones,
including those by Morin, Phillips, Max, and Thurston.
The minimax eversions were automatically generated by
flowing downhill in energy using Brakke's Evolver.
\end{abstract}

\maketitle
\section{A History of Sphere Eversions}

To evert a sphere is to turn it inside-out
by means of a continuous deformation,
which allows the surface to pass through itself,
but forbids puncturing, ripping, creasing, or pinching the surface.
An abstract theorem proved by Smale in the late 1950s implied
that sphere eversions were possible~\cite{Smale}, but it remained a challenge
for many years to exhibit an explicit eversion.
Because the self-intersecting surfaces are complicated
and nonintuitive, communicating an eversion is
yet another challenge, this time in mathematical visualization.
More detailed histories of the problem can be found in~\cite{Levy-MW}
and in Chapter 6 of~\cite{Francis}.

The earliest sphere eversions were designed by hand, and made use
of the idea of a \hwm/.  This is an immersed spherical surface
which is halfway inside-out, in the sense that it has a symmetry interchanging
the two sides of the surface.  If we can find a way to simplify
the \hwm/ to a round sphere, we get an eversion by performing this
simplification first backwards, then forwards again after applying the
symmetry.  The eversions of Arnold Shapiro (see~\cite{FranMor}),
Tony Phillips~\cite{Phillips}, and Bernard Morin~\cite{MorPet}
can all be understood in this way.

\figtwo{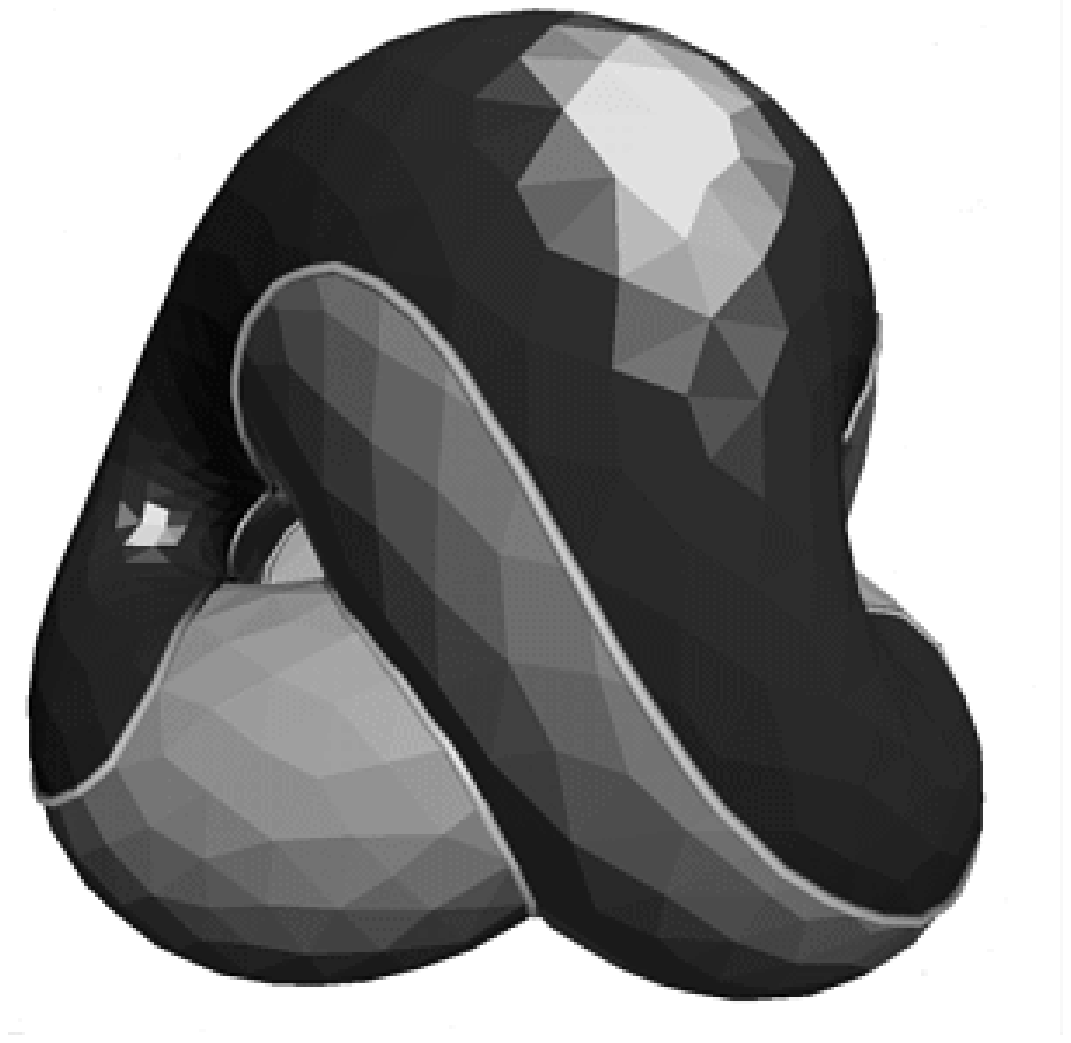}{Halfway models: Boy's surface and Morin surface}
{These are the \hwm/s for the first two minimax eversions.
The Boy's surface (left), an immersed projective plane with three-fold
symmetry and a single triple point, minimizes Willmore's elastic
bending energy.  The figure actually shows an immersed sphere,
double covering Boy's surface, with the two (oppositely-oriented)
sheets pulled apart slightly.  The Morin surface shown (right) also
minimizes Willmore energy; it has a four-fold rotational symmetry
which reverses orientation, exchanging the lighter and darker sides
of the surface.}

In practice, two kinds of \hwm/s have been used, shown in \fig{hwms}.
The first, used by Shapiro and by Phillips (see \fig{phi-max}, left),
is a Boy's surface, which is an immersed projective plane.
In other words it is a way of immersing a sphere in space such that antipodal
points always map to the same place.  Thus there are two opposite sheets
of surface just on top of each other.  If we can succeed in pulling
these sheets apart and simplifying the surface to a round sphere
right-side-out, then pulling them apart the other way will lead to
the inside-out sphere.

The other kind of \hwm/ is a Morin surface; it
has four lobes, two showing the inside and
two the outside.  A ninety-degree rotation interchanges the sides, so
the two halves of the eversion differ by this four-fold twist.
Morin (whose blindness incidentally shows that mathematical
visualization goes well beyond any physical senses)
and Ap{\'e}ry~\cite{AperMor} have shown that an
eversion based on a Morin surface \hwm/ has the minimum
possible number of topological events.

\figtwob{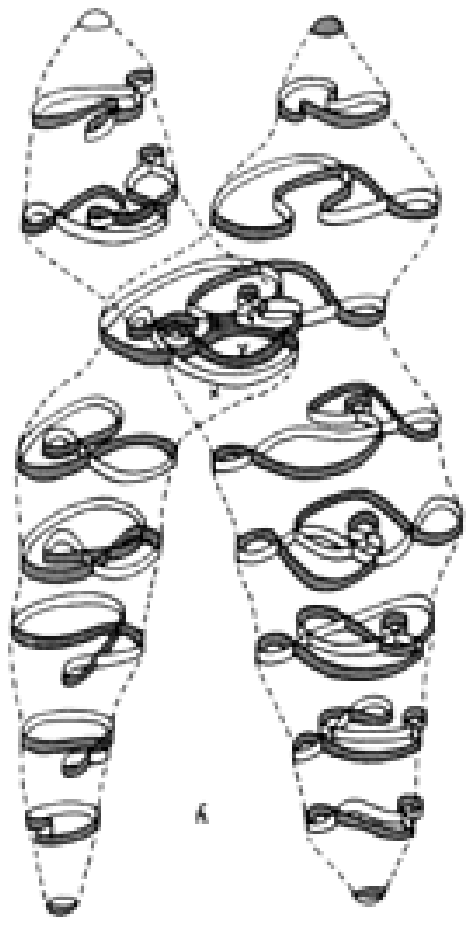}{Phillips drawing and frame from Max video}
{This drawing (left) by Tony Phillips~\cite{Phillips} shows one stage
of his sphere eversion based on a Boy's surface \hwm/.
This frame (right) from Nelson Max's classic computer animation
of a Morin eversion shows a stage near the \hwm/.}

Models of this eversion were made by Charles Pugh, and Nelson Max
digitized these models and interpolated between them for his famous
1977 computer graphics movie ``Turning a Sphere Inside Out''~\cite{max},
a frame of which is shown in \fig{phi-max}.
Morin eversions have also been implemented on computers
by Robert Grzeszczuk and by John Hughes (see \fig{grz-hug}),
among others.  These eversions have usually been laboriously
created by hand, splining or interpolating between key frames,
or looking for algebraic or trigonometric equations.
\figtwo{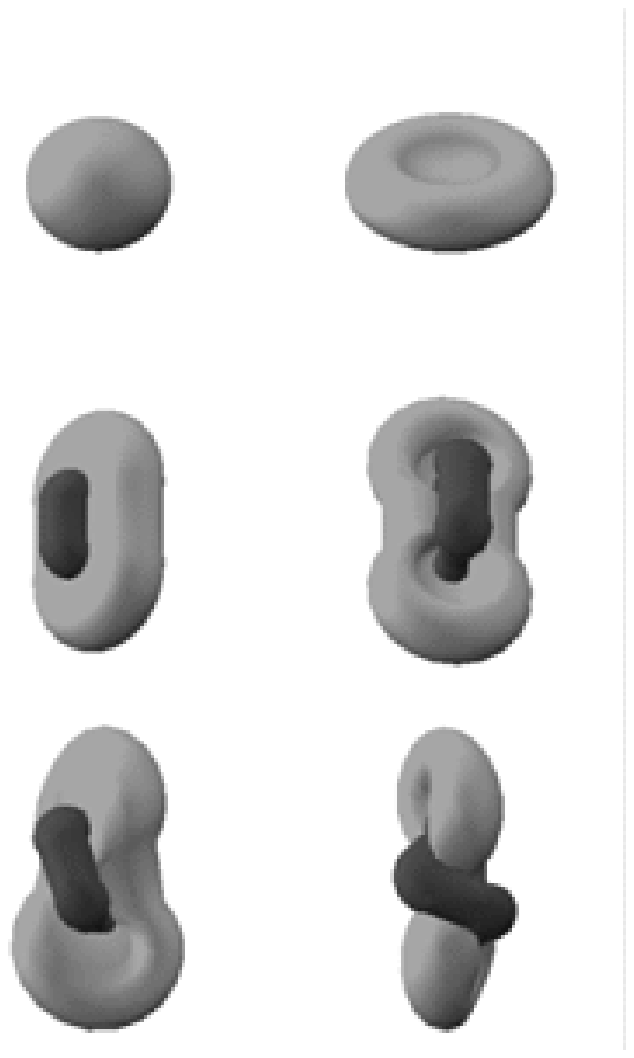}{Scenes from Hughes and Grzeszczuk eversions}
{Still pictures from the Morin eversion implementions by
Robert Grzeszczuk (left) and John Hughes (right).}

A new proof of Smale's original theorem, providing more geometric insight,
was developed by Bill Thurston (see~\cite{Levy-MW}).  His idea was to
take any continuous motion (homotopy) between two surfaces, and make it
regular (taking out any pinching or creasing) when this is possible by adding
corrugations in a controlled way; the corrugations make the surface
more flexible.
This kind of sphere eversion was beautifully illustrated in the
computer-graphics video ``Outside In''~\cite{OutsideIn}
produced at the Geometry Center in 1994 (see \fig{oi}).
\figtwob{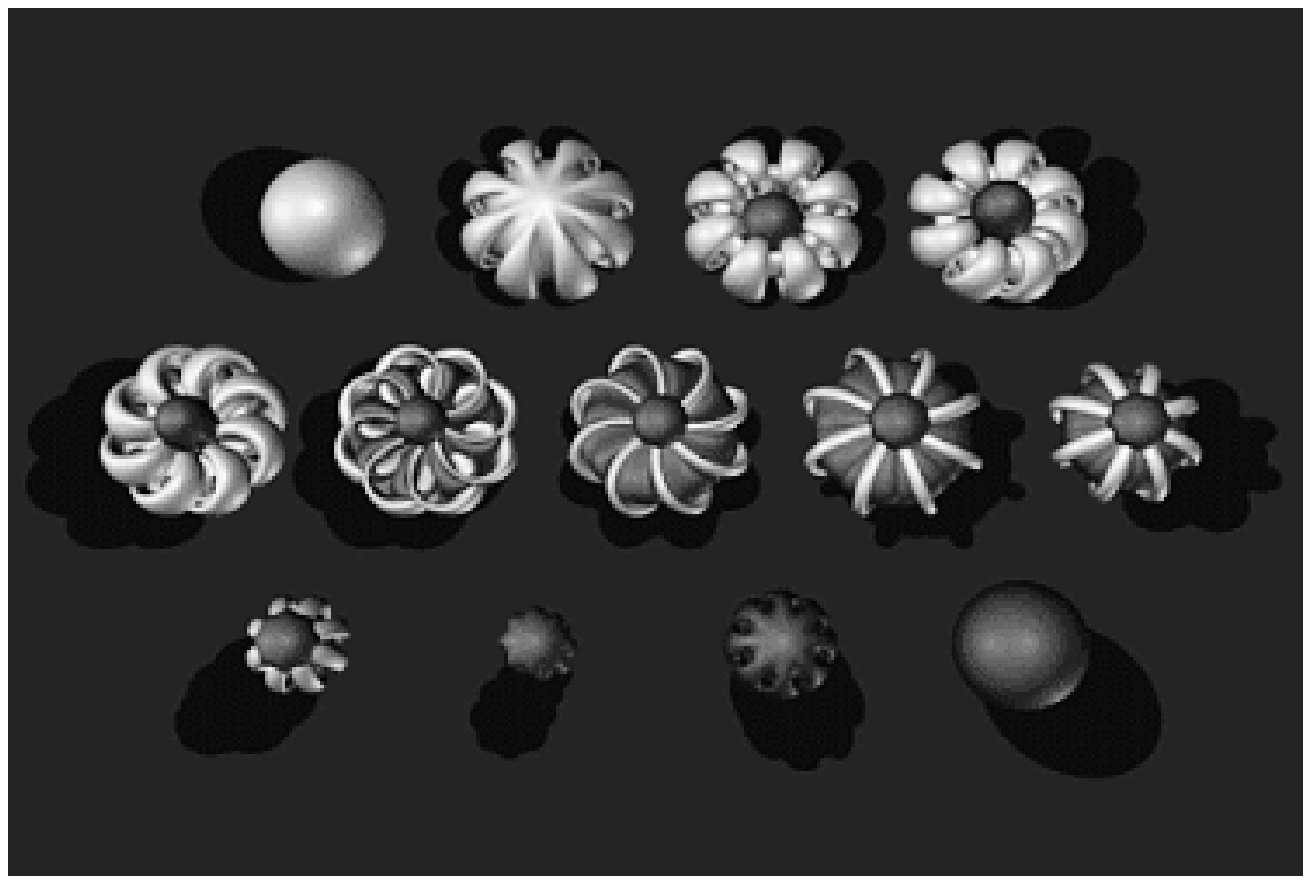}{Scenes from Outside In}
{These pictures, from ``Outside In'', show Thurston's sphere eversion,
implemented through corrugations.}
The corrugation idea is
quite natural, and provides a way to understand all regular homotopies,
not just the sphere eversion.  However, the resulting eversion, though
nicely full of symmetry (except, interestingly, the temporal symmetry
seen in all sphere eversions based on \hwm/s), is quite elaborate,
and has many more topological events than necessary.

\section{Bending Energy and the Minimax Eversions}
Our minimax sphere eversions differ from the earlier ones mentioned
above in that they are computed automatically by a process of energy
minimization.
The elastic bending energy for a stiff wire is the integral
of squared curvature.  For a surface in space, at each point
there are two principal curvatures, and their average, the mean
curvature, shows how much the surface deviates from being
minimal.  The integral of squared mean curvature is thus a bending
energy for surfaces, often called the Willmore
energy~\cite{Willmore-Leuven}.  (Mathematically this is equivalent
to many other formulations because of the Gauss-Bonnet theorem.)

Among all closed surfaces, the round sphere minimizes this bending energy.
(The energy is scale-invariant, and we normalize so the sphere has
energy $1$.)
It is also known~\cite{LiYau,Kus-conf} that any self-intersecting
surface with a $k$-tuple point has energy at least $k$.
To evert a sphere, it is necessary to pass through some stage
with a quadruple point~\cite{BanMax,Hughes-quad}, and hence energy
at least $4$.

Francis and Morin realized that the Morin surface
and Boy's surface are just the first in an infinite sequence of
possible \hwm/s for what they call the tobacco-pouch
eversions (see \fig{gkf-tob1}, taken from~\cite{Francis}).
\figrb{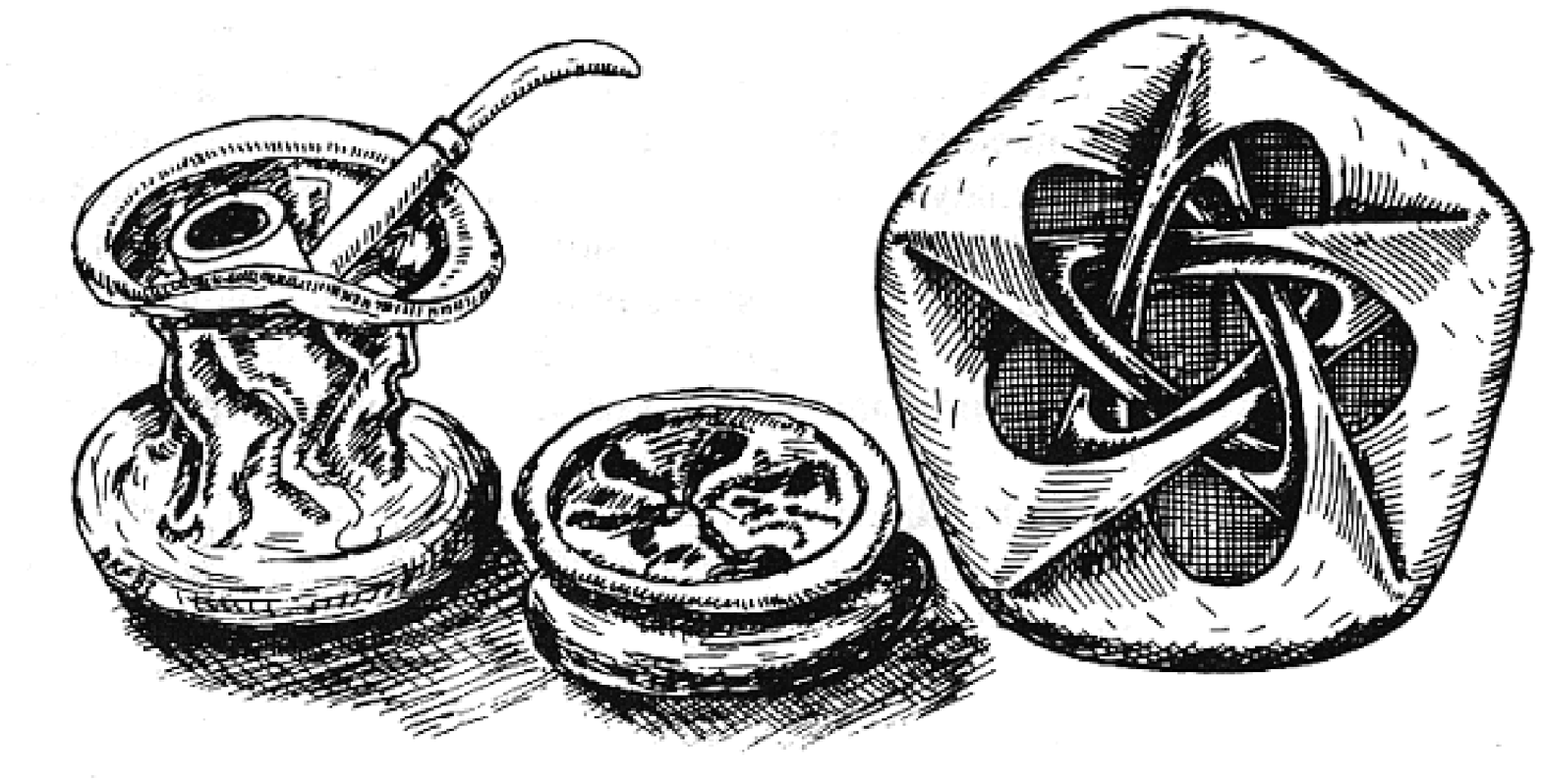}{2.2}{Francis tobacco-pouch sketches}
{These sketches by George Francis show a French tobacco pouch,
and a cutaway view of the \hwm/ for the order-five tobacco-pouch eversion.}
In the 1980's, Robert Bryant~\cite{Bry-dual} classified
all critical points for the Willmore energy among spheres;
they all have integer energy values which (except for the round sphere)
are at least $4$.
Rob Kusner, being familiar with both of these results, realized
that he could find, among Bryant's critical spheres, ones with
the tobacco-pouch symmetries.  Among surfaces with those symmetries,
Kusner's have the least possible Willmore energy.

In particular, Kusner's Morin surface with four-fold
orientation-reversing symmetry has energy exactly~$4$.
If we don't enforce this symmetry, then presumably this surface
is an unstable critical point---a saddle point for the energy.
Pushing off from this saddle in one downhill direction, and flowing down
by gradient descent, we should arrive at the round sphere,
since it is the only critical point with lower energy.  (Of course,
there's not enough theory for fourth-order partial differential equations
for us to know in advance that the surface will remain smooth and not pinch
off somewhere.)

As we saw before,
such a homotopy, when repeated with a twist, will give us a sphere
eversion.  This eversion starts at the round sphere (which minimizes
energy) and goes up over the lowest energy saddle point; then
it comes back down on the other side arriving at the inside-out round sphere.

Around 1995, in collaboration with Kusner and Francis, I computed
this minimax eversion~\cite{FSetal-minimax},
using Ken Brakke's Evolver~\cite{Bra-evol},  which is a software
package designed for solving variational problems, like finding the
shape of soap films or (see~\cite{HKS-willmore})
minimizing Willmore energy.
We were pleased that the computed eversion
was not only geometrically optimal in the sense of requiring the least
bending, but also topologically optimal, in that it was one of
the Morin eversions with the fewest topological events.

Computations of the higher-order minimax sphere eversions~\cite{CSE}
(like the one with a Boy's surface \hwm/) had to wait until Brakke and
I implemented some new symmetry features in the Evolver
(see~\cite{BS-evolsym}).  A minimax tobacco-pouch eversion,
whose \hwm/ has $2k$-fold symmetry, will break some of this symmetry.
But it will maintain $k$-fold symmetry throughout, and the
Evolver now works with only a single fundamental domain for this symmetry.
We can find the initial \hwm/s
by minimizing bending energy while enforcing the full symmetries.

Alternatively, we can compute them directly. Bryant's classification
says that all critical spheres are obtained as conformal transformations
of minimal surfaces, and Kusner gave explicit Weierstrass parameterizations
for the minimal surfaces he needed.  In \fig{minimal} we see the minimal
surface with four flat ends which transforms into our Morin \hwm/. 
\figtwo{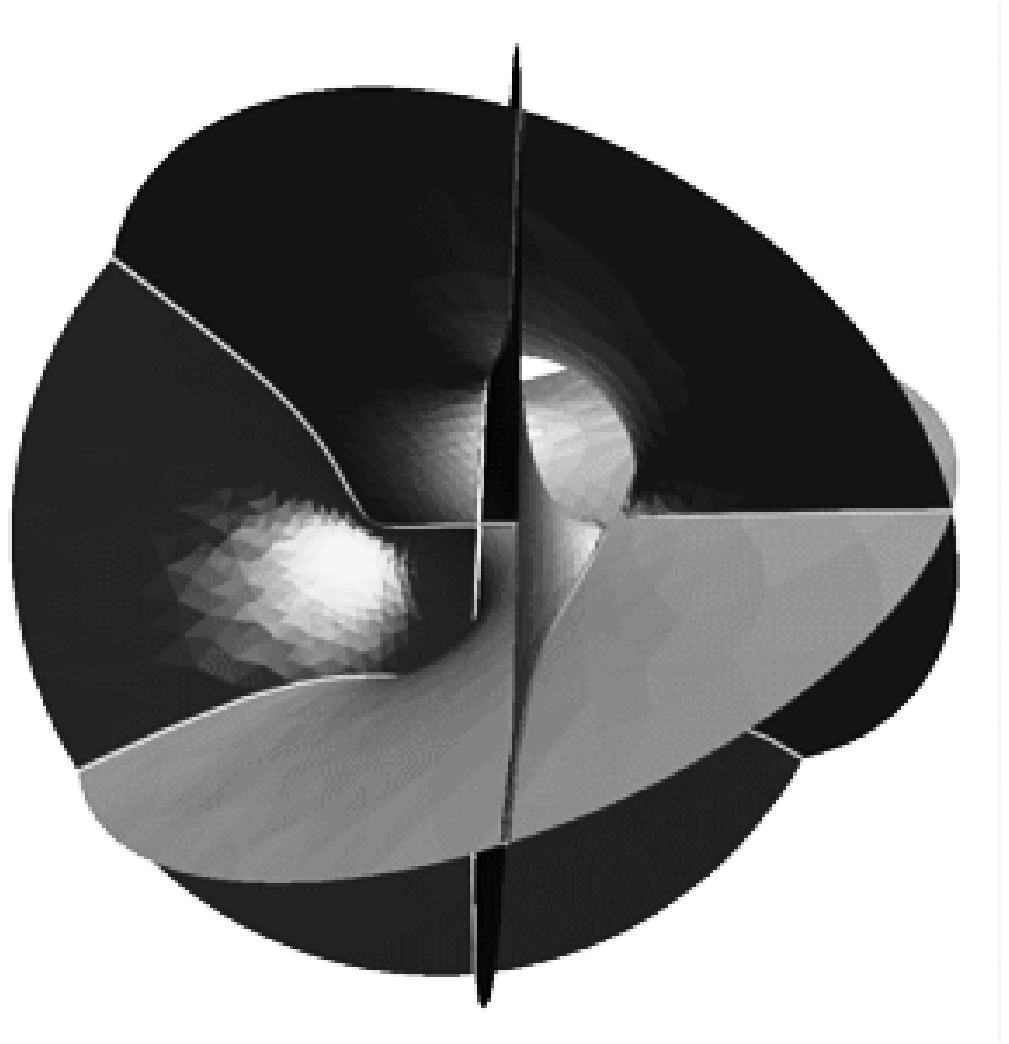}{Kusner's minimal surface}
{This minimal surface (left), with four flat ends, gives rise to Kusner's
Morin surface of least Willmore bending energy, when a conformal
M\"obius transformation is applied to compactify it.  If the transformation
sends the double-tangent point to infinity, we get an interesting
picture (right).}

The Evolver works with triangulated approximations (with a few thousand
triangles) to the true smooth
surfaces, and we update the triangulation as needed to maintain
a good approximation.  Initially, it's necessary to use the second-order
Hessian methods implemented in the Evolver to push off the saddle point
and find our way downhill.

In 1998, working with Francis and Stuart Levy, I produced a computer
graphics video, ``The Optiverse''~\cite{optiverse}, which shows the first
four minimax eversions.\footnote{This video was produced at the NCSA at
the University of Illinois.  More information about it is
available on our website at
http:/\kern-.1em/new.math.uiuc.edu/optiverse/.}
It premiered at the International Congress
of Mathematicians in Berlin, and was also shown at SIGGRAPH and
written up nicely in Science News~\cite{scinews}.

Many scenes in the video (as in \fig{styles}) capture views of
the eversion also available in our real-time interactive computer
animation running on desktop workstations or immersive virtual-reality
environments.
These show a solid, colored surface, with white tubes around
the double-curves of self-intersection.  Here we are aware of the
triangulation used for the computation, especially if we leave gaps
between triangles or in their middles.
\figtwob{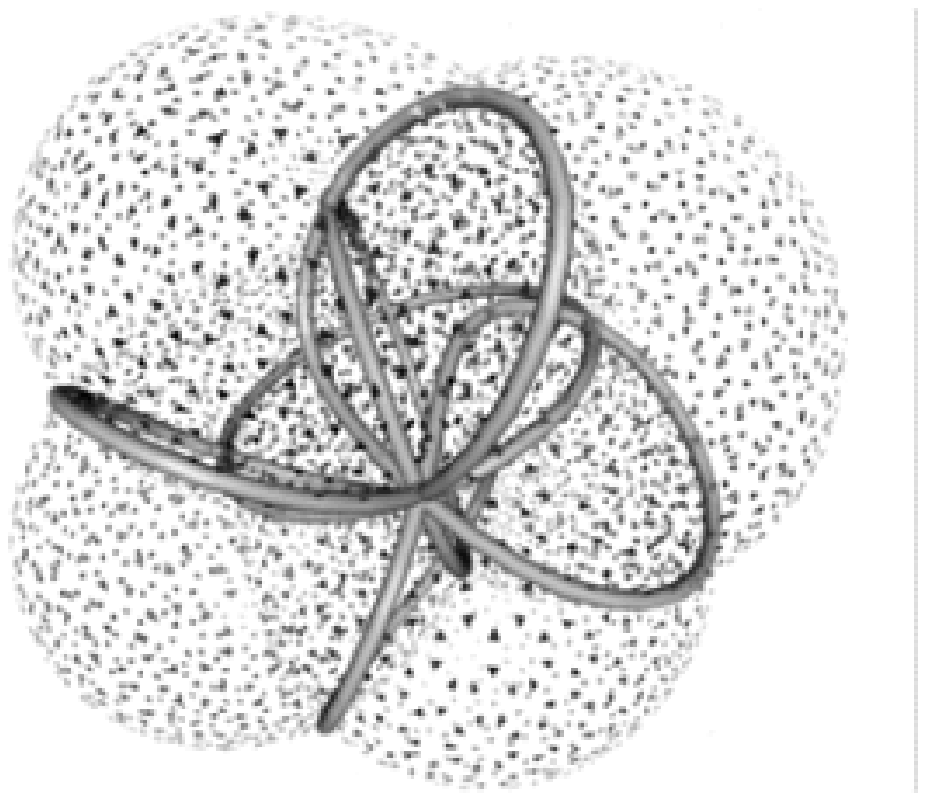}{Gapped halfway model and framework gastrula}
{The Morin halfway model (left), shown with all triangles shrunk, has
an elaborate set of double curves where the surface crosses itself.
An late stage in the eversion is like a gastrula (right), shown here as
a triangular framework.}

But no one method of rendering can give all useful visual information
for a nonintuitive phenomenon like a self-intersecting surface.
So other parts of the video (like \fig{transp})
were rendered as transparent soap films,
with the soap bubble shader I wrote for Renderman (see~\cite{leon}).
\figr{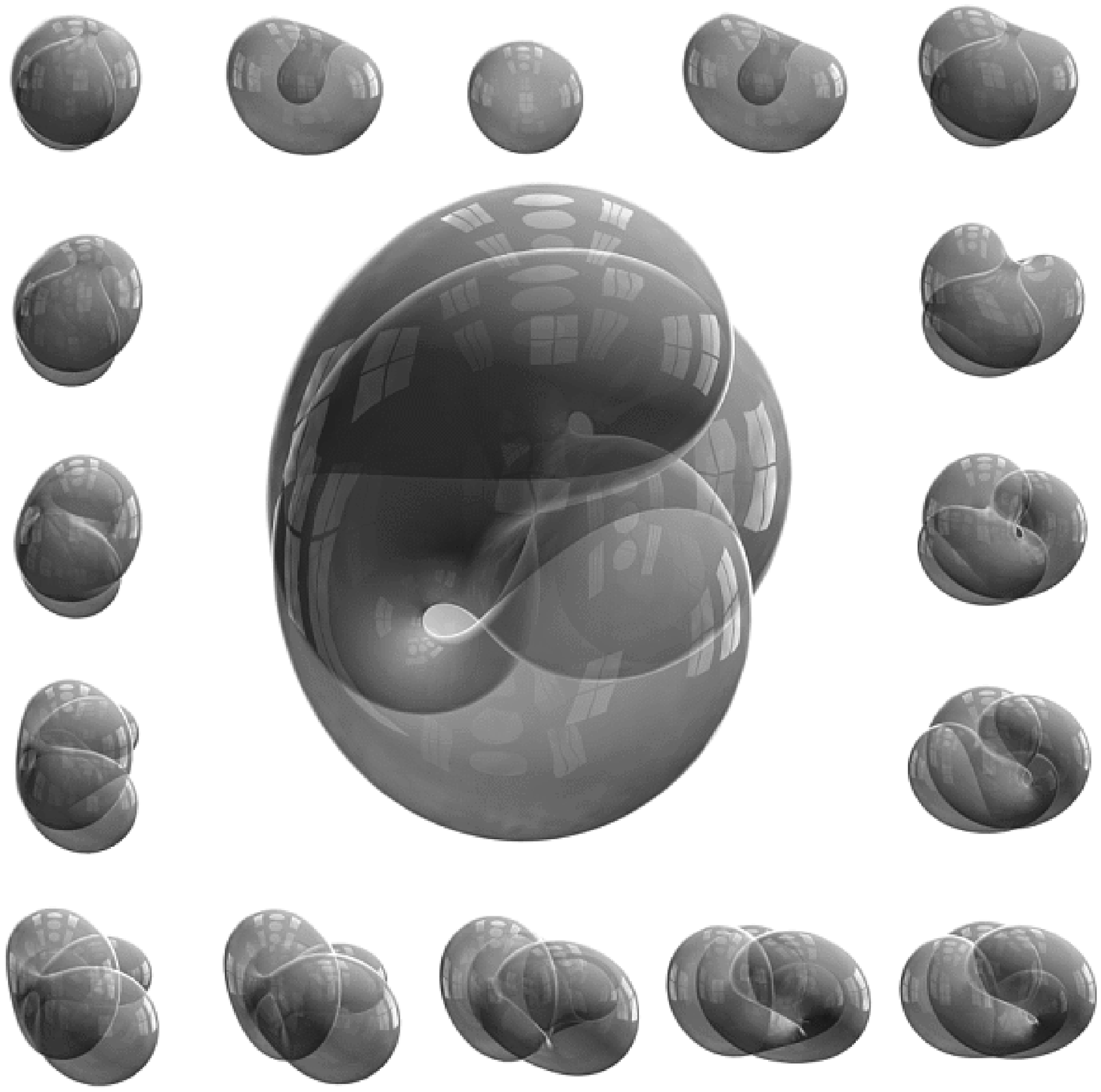}{4.9}{Transparent eversion montage}
{This minimax sphere eversion is a geometrically optimal way
to turn a sphere inside out, minimizing the elastic bending
energy needed in the middle of the eversion.
Starting from the round sphere (top, moving clockwise),
we push the north pole down, then push it through the south pole (upper right)
to create the first double-curve of self-intersection.
Two sides of the neck then bulge up, and these
bulges push through each other (right) to give the second double-curve.
The two double-curves approach each other, and when
they cross (lower right) pairs of triple points are created.
In the \hwm/ (bottom) all four triple points merge at
the quadruple point, and five isthmus events happen simultaneously.
This \hwm/ is a symmetric critical point for the
Willmore bending energy for surfaces.
Its four-fold rotation symmetry interchanges
its inside and outside surfaces.  Therefore, the second half of the eversion
(left) can proceed through exactly the same stages in reverse order, after
making the ninety-degree twist.
The large central image belongs between the two lowest ones on the
right, slightly before the birth of the triple points.}

\section{Topological Stages in the Eversion}

A generic regular homotopy has topological events occuring only
at isolated times, when the combinatorics of the self-intersection
changes.  These events happen at times when the surface normals
at a point of intersection are not linearly independent.  That
means either when there are two sheets of the surface tangent to each other
(and a double curve is created, annihilated, or reconnected),
or when three intersecting sheets share a tangent line (creating or
destroying a pair of triple points), or when four sheets come together
at a quadruple point.

The double tangency events come in three flavors.  These can be modeled
with a rising water level across a fixed landscape.  As the water
rises, we can observe creation of a new lake, the conversion between
a isthmus of land and a channel of water, or the submersion of an
island.  These correspond to the creation, reconnection, or annihilation
of a double-curve, here seen as the shoreline.

\figrp{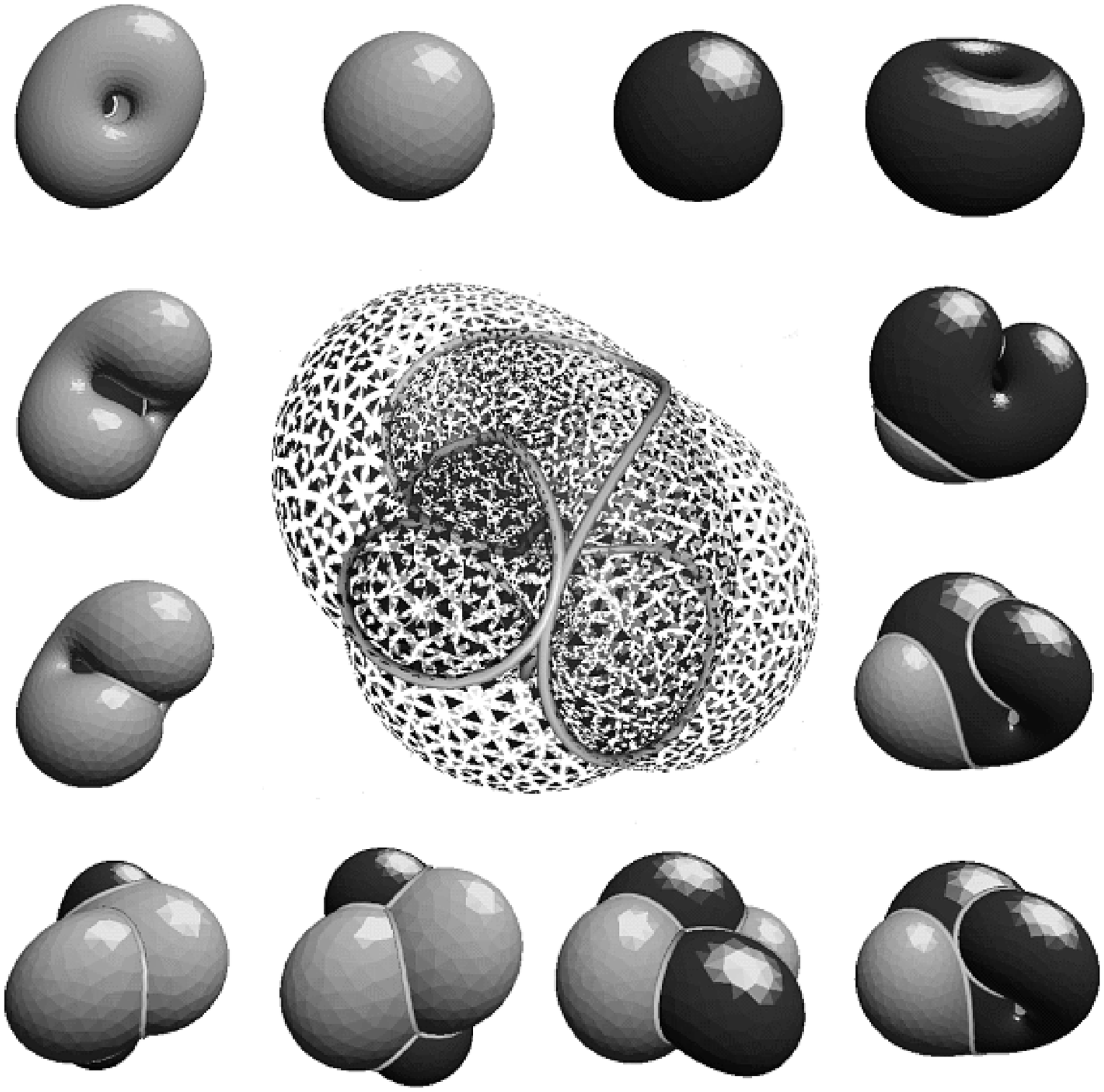}{4.9}{Two-fold eversion montage}
{This is the same eversion shown in \fig{transp}, but rendered with
solid surfaces.  Again, we start at the top with a round sphere, and
proceed clockwise.
Down the right-hand side we see the creation first of two double-curves,
and then of a pair of triple points.  (Another pair is created at
the same time in back; the eversion always has two-fold rotational
symmetry.) Across the bottom, we go through the \hwm/,
interchanging the roles of the dark and light sides of the surface.
Up the left column, we see the double-curves disappear one after
the other.  (This time, the figures on the left are not exactly the
same stages as the corresponding ones on the right.)
In the center, we examine the double locus just when
pairs of triple points are being created, by shrinking each triangle
of the surface to a quarter of its normal size.}

In the minimax eversion with two-fold symmetry,
seen in \fig{transp} and \fig{bivtile},
the first two events create the two double-curves.  When these
twist around to intersect each other, two pairs of
triple points are created.  At the halfway stage, six events happen
all at the same time.  Along the symmetry axis, at one end we have
a quadruple point, while at the other end we have a double tangency
creating an isthmus event.  Finally, at the four ``ears'', at the
inside edge of the large lobes, we have additional isthmus events:
two ears open as the other two close.  (See~\cite{FSetal-minimax}
for more details on these events.)

The three-fold minimax eversion, using the Boy's surface of
\fig{hwms} as a \hwm/, has too many topological events to describe
easily one-by-one, and its three-fold symmetry means that the
events are no longer all generic.  But we can still follow
the basic outline of the eversion from \fig{tritile}.
\figrp{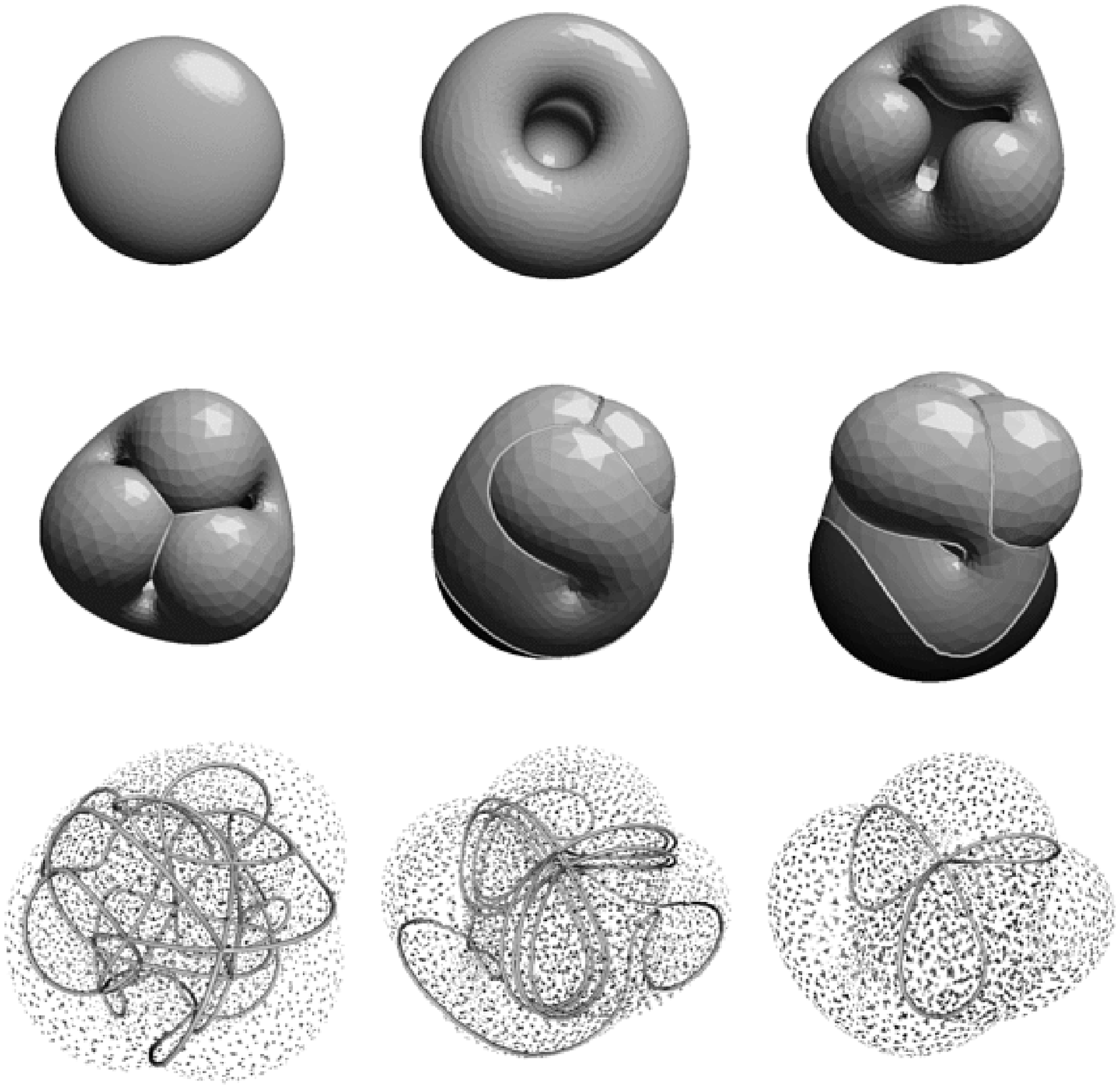}{4.9}{Three-fold eversion}
{This three-fold minimax eversion starts (top row, left-to-right)
with a gastrula stage like that of the two-fold eversion, but
the three-fold symmetry means that three fingers reach up from the neck
instead of two.  They intersect each other (middle row, left-to-right)
and then twist around, while complicated things are happening inside.
The bottom left image is at the same stage as the middle right image,
but with triangles shrunk to show the elaborate double-curves.
These double-curves separate into two pieces (bottom row, left-to-right),
one of which is a four-fold cover of the propeller-shaped double-curve
seen in the Boy's surface \hwm/ (lower right).}

\section*{Acknowledgments}
The minimax sphere eversions described here are work
done jointly in collaboration with
Rob Kusner, Ken Brakke, George Francis, and Stuart Levy,
to whom I owe a great debt.  I would also like to thank
Francis, Robert Grzeszczuk, John Hughes, AK Peters, Nelson Max
and Tony Phillips for permission to reproduce figures from earlier
sphere eversions.  This paper will appear in the proceedings of
two summer 1999 conferences on Mathematics and Art:
ISAMA 99 (June, San Sebasti\'an, Spain), and Bridges (July, Winfield, Kansas).
I wish to thank Nat Friedman, Reza Sarhangi, and Carlo
S\'equin for the invitations to speak at these conferences.
My research is supported in part by NSF grant DMS-97-27859.
\vfill\clearpage
\small
\bibliographystyle{math}
\bibliography{vismath}

\end{document}